# PERCOLATION ON FINITE GRAPHS AND ISOPERIMETRIC INEQUALITIES

By Noga Alon[1], Itai Benjamini and Alan Stacey

*Tel Aviv University, Microsof Research and Weizmann Institute, and Centre for Mathematical Sciences*

Consider a uniform expanders family $G_n$ with a uniform bound on the degrees. It is shown that for any $p$ and $c > 0$, a random subgraph of $G_n$ obtained by retaining each edge, randomly and independently, with probability $p$, will have at most one cluster of size at least $c|G_n|$, with probability going to one, uniformly in $p$. The method from Ajtai, Komlós and Szemerédi [*Combinatorica* **2** (1982) 1–7] is applied to obtain some new results about the critical probability for the emergence of a giant component in random subgraphs of finite regular expanding graphs of high girth, as well as a simple proof of a result of Kesten about the critical probability for bond percolation in high dimensions. Several problems and conjectures regarding percolation on finite transitive graphs are presented.

**1. Introduction.** In this paper we primarily consider percolation on finite graphs and, in particular, the existence and uniqueness of large components, typically meaning components whose size is proportional to the number of vertices in the graph. Our main results in this context apply to expanders, which are graphs satisfying a particular isoperimetric inequality, although we conjecture that these results hold somewhat more generally. The techniques we use can also be used to give a significantly shorter proof than those previously known for the fact that the critical probability for percolation on $\mathbb{Z}^d$ is asymptotically $1/(2d)$ as $d \to \infty$.

Given a graph $G$, we shall use $G(p)$ to denote the spanning subgraph of $G$ obtained by retaining each edge of $G$ independently with probability $p$. This has been very extensively studied in the case when $G$ is a complete graph, and this is known as the standard random graph model or the mean

Received October 2002; revised June 2003.
[1]Supported in part by a USA–Israeli BSF grant, by the Israel Science Foundation and by the Hermann Minkowski Minerva Center for Geometry at Tel Aviv University.

field model; see, for example, the books [9] and [17]. Percolation on general infinite graphs has been studied (see [7] or [21] for background) and there, as in this paper, isoperimetric inequalities play a key role. Most other studies of percolation on finite graphs concern specific graphs, such as the torus, which are closely related to percolation on corresponding infinite graphs such as $\mathbb{Z}^d$. Another example of this phenomenon is the study of the contact process on finite trees [25]; the contact process on a graph $G$ is loosely analogous to percolation on the Cartesian product $G \times \mathbb{Z}$; and both the contact process on $\mathbb{T}$ and percolation on $\mathbb{T} \times \mathbb{Z}$, where $\mathbb{T}$ is a homogeneous tree, have been widely studied.

In light of the above it is, perhaps, surprising that there has been little work regarding percolation on general classes of finite graphs. In this paper we hope to demonstrate that there are interesting questions in this area. The questions asked and methods used draw on the theories of both random (finite) graphs and percolation on infinite graphs.

In two widely studied cases, where $G$ is either the complete graph or resembles a finite subset of $\mathbb{Z}^d$ (either a large $d$-dimensional $n \times \cdots \times n$ torus or box), uniqueness results for the giant component are known. Very precise results are known for the complete graph (see [17] for a recent account). For the torus or box, results can be deduced from information about the corresponding infinite graph; see, for example, Lemma 2 of [11]. It seems natural to conjecture that this uniqueness is a much more general phenomenon.

CONJECTURE 1.1. Let $G_n = (V_n, E_n)$ be a sequence of connected finite transitive graphs with a uniformly bounded maximum degree and with $|V_n| \nearrow \infty$. Suppose that $\mathrm{diameter}(G_n) = o(|V_n|/\log|V_n|)$. Then for any $a > 0$,

$$\sup_p \mathbb{P}_p(\text{there is more than one connected component of size at least } a|G_n|) \to 0$$

as $n \to \infty$, where $\mathbb{P}_p$ denotes the probability with respect to the measure $G(p)$.

It is easily seen, by considering cycles or the Cartesian product of a large cycle with a small transitive graph, that the conjecture fails with the condition on the diameter dropped. These examples fail only for values of $p$ approaching 1, and slightly more sophisticated examples show that it fails without the diameter condition even with $p$ bounded away from 1. Indeed, the product of a regular expander of order $c \log n$ and a cycle of order $n/(c \log n)$ forms such a family of examples. Similarly, the product of a complete graph and a triangle shows that the assumption on the bounded degrees is also essential.

Our first result, in Section 2, establishes uniqueness of the giant component for expanders. This holds even without vertex transitivity since the expansion property gives sufficient uniform control over the geometry of the graph. In fact, slightly more can be shown: the uniqueness holds in this case even for clusters of sublinear size; see Theorem 2.8 for a detailed statement.

In [1], Ajtai, Komlós and Szemerédi proved that the critical probability for the emergence of the giant component in bond percolation on the hypercube $\{0,1\}^d$ is $1/d$. The strategy of the proof is twofold. First, one uses the very local geometry of the hypercube in a neighborhood of a vertex to obtain (based on a basic branching process argument) that percolation with $p > 1/d$ will have many clusters of polylogarithmic size. These clusters cover a constant fraction of the hypercube. In the second step one uses the isoperimetric inequality for the hypercube to prove that by adding additional independent $\varepsilon/d$ percolation, most of these polylogarithmic clusters join to form a giant component. In Section 3 we remark on how this approach can be naturally used to determine the critical probability for percolation on some other graphs, including regular expanders with large girth. This technique also enables us to present, in Section 4, a rather simple proof for the fact [18] that the critical probability for bond percolation in $\mathbb{Z}^d$ is $\frac{1+o(1)}{2d}$, as $d \to \infty$.

1.1. *Expanders and other definitions.* Expanders are defined in terms of a certain isoperimetric inequality. Such inequalities have wide applications in graph theory and in percolation in particular. They play a crucial role in the study of percolation on general infinite graphs; for a few natural conjectures relating isoperimetric inequalities to percolation in this context see [7], especially Conjecture 1, Question 2. Although considerable progress has been made in recent years, there is scope for further work in understanding the relation between properties of percolation processes and the isoperimetric profile of the underlying graph (in the spirit in which the behavior of the simple random walk is directly linked to isoperimetric inequalities, see, e.g., [13]).

Another important example of the use of isoperimetric inequalities in the area of graph theory is the role played by *conductance*, an isoperimetrically-defined quantity, in showing that Markov chains are rapidly mixing. See, for example, [6].

We now turn to the precise definitions. The *girth* $g(G)$ of a graph $G = (V, E)$ is the minimum length of a cycle in $G$. For any two sets of vertices in $G$, $A, B \subseteq V$, the set $E(A, B)$ consists of all those edges with one endpoint in $A$ and the other in $B$. For a finite graph $G$ its edge-isoperimetric number, $c(G)$, (also called its Cheeger constant) is given by

$$\min_{\substack{A \subset V \\ 0 < |A| \leq |V|/2}} \frac{|E(A, V \setminus A)|}{|A|}.$$

We will also make use of the vertex isoperimetric constant, $\iota(G)$, which we now define similarly. Given a set of vertices $A \subseteq V$, define the external boundary of $A$, $\partial A$, to consist of those vertices outside $A$ which have a neighbor in $A$. Then define

$$\iota(G) = \min_{\substack{A \subset V \\ 0 < |A| \leq |V|/2}} \frac{|\partial A|}{|A|}.$$

We shall be interested in families of graphs whose isoperimetric constants are bounded away from 0. Given $b > 0$, we say that a graph, $G$, is an *edge b-expander* if $c(G) \geq b$ and a *vertex b-expander* if $\iota(G) \geq b$. We shall also refer, with a slight abuse of notation, to a set of graphs, or a sequence of graphs $(G_n)$, as an *edge* (resp. *vertex*) *b-expander* if each graph in the set is an edge (resp. vertex) $b$-expander. A sequence of graphs is called simply an *edge* (resp. *vertex*) *expander* if it is an edge (resp. vertex) $b$-expander for some $b > 0$. Most sequences we consider will have a uniform bound, $\Delta$, say, on the degrees of the vertices, and in that case it is clear that the sequence is a vertex expander if and only if it is an edge expander; we refer to such sequences simply as *expanders*.

Expanders received a considerable amount of attention in the literature in recent years, mostly because these graphs have numerous applications in theoretical computer science; see, for example, [4, 20]. It is well known that for any fixed $d > 2$, random $d$-regular graphs of size $n$ are asymptotically almost surely expanders, as $n$ grows. The problem of constructing infinite families of bounded degree expanders is more difficult, and there are several known constructions of this type. Most of these constructions are Cayley graphs, and are therefore vertex transitive.

The distance between two vertices of a graph is the length of the shortest path between them. Given a vertex $v$, the set of vertices within distance $r$ from $v$ (or the subgraph they induce) will be denoted by $B(v, r)$. Also, for a set of vertices $A$, $B(A, r)$ will denote the set of all vertices which are within distance $r$ of some vertex in $A$.

**2. Uniqueness of the giant component.** The aim of this section is to establish Conjecture 1.1 with the condition of vertex transitivity replaced by the condition of expansion.

THEOREM 2.1. *Let $b > 0$ and let $\Delta \in \mathbb{N}$. Let $G_n = (V_n, E_n)$ be a sequence of graphs with maximum degree at most $\Delta$ which are vertex b-expanders, with $|V_n| \to \infty$. Let $0 \leq p_n \leq 1$ and let $c > 0$. Then*

$$\mathbb{P}(G_n(p_n) \text{ contains more than one component of order at least } c|V_n|) \to 0 \quad (1)$$

*as $n \to \infty$.*

The statement of this theorem holds for any family of expanders, such as the ones described in [4, 20, 24] and their references. Various applications of expanders rely on their fault-tolerance as networks that imply that even after deleting an appropriate constant fraction of their edges (arbitrarily), the remaining graphs still contain some linear size connected components or some linear size paths; see, for example, [3, 26]. The theorem above provides more information in the case when the edges are deleted by a random process.

We will refer to components of order at least $c|V_n|$ as *large*. Note that if $p_n \leq a$ for some $a < 1/\Delta$, then standard branching process arguments (see, e.g., [16]) show that the probability of the existence of any large component tends to zero as $n \to \infty$. We use the following lemma to deal with the case when $p$ is close to 1.

LEMMA 2.2. *Let $b > 0$ and let $G_n = (V_n, E_n)$ be a sequence of graphs with maximum degree at most $\Delta$ which are edge $b$-expanders, with $|V_n| \to \infty$. Let $A > 0$ be such that $(\Delta e)A^b < 1$ and let $1 - A \leq p_n \leq 1$ for each $n$. Then for any $c > 0$,*

(2) $\quad \mathbb{P}(G_n(p_n) \text{ contains a component of size between } c|V_n| \text{ and } (1/2)|V_n|) \to 0$

*as $n \to \infty$.*

PROOF. There are various ways to see this, but simple counting turns out to be the most useful in what follows. We shall say that a subset of the vertices of a graph $G$ is *connected* (*in $G$*) if the subgraph of $G$ it induces is connected. We shall make use of the fact (see, e.g., [2]) that in a graph $G = (V, E)$ of maximum degree $\Delta$, the number of connected subsets of $V$ of size $r$, containing some given vertex, is at most $(\Delta e)^r$. Summing over all the vertices counts each such subset $r$ times, so the total number of connected (in $G_n$) subsets of $V_n$ of order $r$ is at most

$$\frac{|V_n|}{r}(\Delta e)^r.$$

Now for any subset, $U$, of $V_n$ of size $r$, where $r \leq |V_n|/2$, the expansion property gives that $|E(U, U^c)| \geq br$, so the probability that all the edges of $E(U, U^c)$ are absent from $G_n(p_n)$ is at most $(1 - p_n)^{br}$; this is an upper bound on the probability that $U$ is the vertex set of some connected component of $G_n(p_n)$. Therefore, the probability there is a component of size between $c|V_n|$ and $|V_n|/2$ is at most

(3) $$\sum_{r=\lceil c|V_n|\rceil}^{\lfloor |V_n|/2\rfloor} \frac{|V_n|}{r}(\Delta e)^r (1-p_n)^{br} \leq \frac{1}{c}\frac{((\Delta e)A^b)^{c|V_n|}}{1 - (\Delta e)A^b},$$

using the fact that $(\Delta e)A^b < 1$. The upper bound in (3) clearly tends to 0 as $|V_n| \to \infty$ and establishes the lemma. □

Turning back to Theorem 2.1, the vertex $b$-expanders in the statement of the theorem must also be edge $b$-expanders (with the same value of $b$). Taking $A$ as in Lemma 2.2, and using the fact that if there are two components of size at least $c|V_n|$, one of them must contain no more than half the vertices, we see that if $p_n \geq 1 - A$, then the probability there is more than one large component is small when $n$ is large. Combining this with the observation about small $p_n$ made after the statement of Theorem 2.1, we may assume that $p_n \in [x, 1-x]$ for all $n$, where $x = \min(1/(2\Delta), 1 - A)$.

In the following very useful lemma (which we do not claim is new), recall that a subset, $\mathcal{X}$, of $\mathcal{P}(E)$ is an *up-set* if, whenever $A \in \mathcal{X}$ and $A \subset B \subseteq E$, then $B \in \mathcal{X}$.

LEMMA 2.3. *Let $x > 0$. Then there exists $\alpha > 0$ so that the following holds. Let $E$ be a finite set and let $\mathcal{A} \subseteq \mathcal{P}(E)$ be an up-set. Given $A \subseteq E$ and $e \in E$, say that $e \in E$ is $A$-pivotal (for $\mathcal{A}$) if $A^e = A \cup \{e\} \in \mathcal{A}$ and $A_e = A \setminus \{e\} \notin \mathcal{A}$. Let $A \subseteq E$ be obtained by selecting each element of $E$ independently with probability $p$, where $p \in [x, 1-x]$, and let $e \in E$ be chosen uniformly at random and independently of the choice of $A$. Then*

$$(4) \qquad \mathbb{P}(e \text{ is } A\text{-pivotal for } \mathcal{A}) \leq \frac{\alpha}{\sqrt{|E|}}.$$

PROOF. Given $E$ and $p \in [x, 1-x]$, we construct a pair $(A, e)$ as follows. Order the elements of $E$ randomly, $e_1 < e_2 < \cdots < e_k$, with each of the $k!$ possible orderings equally likely (with $k = |E|$). Let $X \sim B(k, p)$ be a binomial random variable, independent of this ordering, and let $A$ be the first $X$ elements in the ordering, $A = \{e_1, \ldots, e_X\}$. Now, with probability $X/k$, let $e = e_X$ and with probability $(k - X)/k$, let $e = e_{X+1}$. We now see why this construction yields a pair $(A, e)$ with the distribution given in the statement of the lemma. The fact that the marginal distribution of $A$ is correct is immediate. Now, for any $A$ that arises in this way, it is equally likely to have arisen from any of the $|A|!(k - |A|)!$ orderings which place the elements of $A$ in the first $|A|$ places; the proportion of these orderings in which any given element of $A$ is in the $|A|$th place is exactly $1/|A|$, so the chance that any given element of $A$ turns out to be the random element $e$ is exactly $(1/|A|)(|A|/k)$ or $1/k$. Similarly, again conditional on the choice of $A$, any element outside $A$ also has chance $1/k$ to be equal to $e$.

This apparently peculiar way of constructing $(A, e)$ is useful in estimating the probability that $e$ is $A$-pivotal. Having chosen the ordering, let $A_l = \{e_1, \ldots, e_l\}$, for $0 \leq l \leq k$. Since $\mathcal{A}$ is an up-set there will be precisely one $l$ with the property that $A_l \notin \mathcal{A}$ and $A_{l+1} \in \mathcal{A}$. We see that $e$ is $A$-pivotal precisely if $e = e_{l+1}$ and $X = l$ or $X = l+1$. This happens with probability (conditional on the ordering)

$$\frac{k-l}{k}\mathbb{P}_p(X=l) + \frac{l+1}{k}\mathbb{P}_p(X=l+1) \leq \frac{k+1}{k}\max_{p,m}\mathbb{P}_p(X=m),$$

where the maximum is taken over all $p \in [x, 1-x]$ and all $m$. Since this bound is independent of the ordering, it is also a bound on the unconditional probability than $e$ is $A$-pivotal. However, it is well known (it follows fairly easily, e.g., from bounds on binomial coefficients given by (1.5) in [9]) that this maximum is bounded above by a constant over $\sqrt{k}$, the precise constant depending on the value of $x$. □

Given a subgraph, $H$, of $G_n$, we say that an edge $e \in E(G_n)$ is an $L$-*bridge* if $H_e$ contains two large components which are connected by $e$. Recall that the definition of large depends on the choice of some constant $c > 0$.

COROLLARY 2.4. *Let $x > 0$, $b > 0$, $c > 0$ and $\Delta \in \mathbb{N}$ be given. Then there exists $\beta > 0$ so that the following holds. Let $G_n$ be a graph satisfying the conditions of Theorem* 2.1, *and let $p_n \in [x, 1-x]$. For $e \in E_n$ let $S(e,c,n)$ be the event that $e$ is an $L$-bridge in $G_n(p_n)$. Let $S(c,n)$ be the event that $S(e_n, c, n)$ occurs for an edge $e_n$ chosen uniformly at random. Then*

(5) $$\mathbb{P}(S(c,n)) \leq \beta/\sqrt{|V_n|}.$$

PROOF. We could proceed by adapting the proof of Lemma 2.3, noting that given any ordering on the edges as in that proof, at most $\lfloor 1/c \rfloor$ edges, $e_l$, can be $L$-bridges for the corresponding configuration $A_l$ or $A_{l+1}$. However, since Lemma 2.3 is an attractive general result we prefer to proceed by applying this lemma directly. We do this by constructing up-sets in such a way that any $L$-bridge is pivotal for one of these up-sets.

Given any configuration of edges, $F \subseteq E_n$, let $Y(F)$ count the number of vertices which belong to large components, and let $C(F)$ count the number of large components. Now set

$$Z(F) = \frac{Y(F)}{c|V_n|} - C(F).$$

It is not hard to see that $Z$ is an increasing function of $F$: all we have to note is that if the addition of an edge increases $C(F)$ by 1, then at the same time $Y(F)$ must increase by at least $c|V_n|$. Therefore, for any $t$, the set of configurations, $F$, satisfying $Z(F) \geq t$ is an up-set.

Now let $\mathcal{A}_i = \{F : Z(F) \geq i\}$ for $i = 1, 2, \ldots, \lfloor 1/c \rfloor - 1$ [noting that the maximum value of $Z(F)$ is $1/c - 1$]. Then any $L$-*bridge* is pivotal for some $\mathcal{A}_i$. Hence, applying Lemma 2.3 and summing over $i$,

$$\mathbb{P}(S(c,n)) \leq (\lfloor 1/c \rfloor - 1) \frac{\alpha}{\sqrt{|E_n|}}; \tag{6}$$

but since any expander is connected, $|E_n| \geq |V_n| - 1$, giving (5). □

COROLLARY 2.5. *Let $x > 0$, $b > 0$, $c > 0$, $r > 0$ and $\Delta \in \mathbb{N}$ be given. Then there exists $\gamma > 0$ so that the following holds. Let $G_n$ be a graph satisfying the conditions of Theorem 2.1, and let $p_n \in [x, 1-x]$. Let $w_n$ be a vertex chosen uniformly at random from $V_n$ and let $S'(c,n,r)$ be the event that there is an edge $e$, contained in the ball $B(w_n, r)$, for which the event $S(e,c,n)$ occurs. Then*

$$\mathbb{P}(S'(c,n,r)) \leq \gamma/\sqrt{|V_n|}. \tag{7}$$

PROOF.

$$\mathbb{P}(S'(c,n,r)) \leq \sum_e \mathbb{P}(e \in B(w_n, r)) \mathbb{P}(S(e,c,n))$$

$$\leq \max_e \mathbb{P}(e \in B(w_n, r)) \sum_e \mathbb{P}(S(e,c,n))$$

$$= \max_e \mathbb{P}(e \in B(w_n, r)) |E_n| \mathbb{P}(S(c,n))$$

$$\leq \frac{\Delta^r}{|V_n|} |E_n| \mathbb{P}(S(c,n)).$$

However, $|E_n| \leq |V_n| \Delta/2$, so applying Corollary 2.4 we see that (7) holds with $\gamma = \Delta^{r+1} \beta/2$. □

We now work towards establishing a lower bound for the probability in (7) in terms of the probability of the existence of two or more large components.

LEMMA 2.6. *Given $b > 0$, $c > 0$ and $k < 1$, there exists $r \in \mathbb{N}$ with the following property. Let $G = (V, E)$ be a vertex $b$-expander with $|V| = n$. Suppose that $A \subseteq V$ with $|A| \geq cn$. Then $|B(A, r)| \geq kn$.*

PROOF. For a given $r$, let $C = V \setminus B(A, r)$, and suppose that $|B(A, r)| < kn$ so $|C| > (1-k)n$. Without loss of generality suppose that $c < 1/2$ and $k > 1/2$. By the expansion property, $|B(C, \lceil \log_{1+b}(1/2(1-k)) \rceil)| > n/2$ and $|B(A, \lceil \log_{1+b}(1/2c) \rceil)| \geq n/2$ so these two balls have a vertex in common. Therefore, $A$ and $C$ are within distance $\lceil -\log_{1+b}(2(1-k)) \rceil + \lceil -\log_{1+b}(2c) \rceil$, giving a contradiction if $r$ is greater than this value. □

LEMMA 2.7. *Let $x > 0$, $b > 0$, $c > 0$ and $\Delta \in \mathbb{N}$, and take $r$ as in Lemma 2.6 corresponding to the case $k = 3/4$. Let $G_n$ be as in the statement of Theorem 2.1 and let $x \leq p_n \leq 1 - x$. Let*

$$\delta_n = \mathbb{P}(G_n(p_n) \text{ contains more than one large component}).$$

*Let $S'(c, n, r)$ be as in the statement of Corollary 2.5. Then*

(8) $$\mathbb{P}(S'(c, n, r)) \geq \tfrac{1}{2} x^{2r} \Delta^{-2r^2} \delta_n.$$

PROOF. Let $G_n = (V_n, E_n)$, so bond percolation on $G_n$ simply assigns probabilities to subsets of $E_n$, which we refer to as *configurations*. For each $w \in V_n$, let $D(w, r)$ be the event that the ball $B(w, r)$ contains two vertices belonging to different large components; we shall regard $D(w, r)$ as a subset of the configurations. We shall let $D(r)$ be the event that $D(w, r)$ occurs for a vertex $w$ chosen uniformly at random, so

(9) $$\mathbb{P}(D(r)) = \frac{1}{|V_n|} \sum_{w \in V_n} \mathbb{P}(D(w, r)).$$

Let $S'(w, c, n, r)$ be the event that $B(w, r)$ contains an $L$-bridge for the graph $G_n(p)$. So $S'(c, n, r)$ is the event that $S'(w, c, n, r)$ occurs for a vertex $w$ chosen at random. For each $w$, a configuration in $D(w, r)$ can be transformed into a configuration lying in $S'(w, c, n, r)$ by the addition of some edges lying within $B(w, r)$: indeed, if $D(w, r)$ occurs, then pick two large components which have vertices lying in $B(w, r)$ and choose a shortest path between them; add edges from this path to the configuration until there is a path of edges in the configuration between the two components; the last edge added is then an $L$-bridge. This procedure gives us a function, $f$ say, from $D(w, r)$ (regarded as a set of configurations) to $S'(w, c, n, r)$. Since the function adds at most $2r$ edges, all taken from $B(w, r)$ which contains at most $\Delta^r$ edges, the inverse image of any set in $S'(w, c, n, r)$ contains at most $\binom{\Delta^r}{2r} + \binom{\Delta^r}{2r-1} + \cdots + \binom{\Delta^r}{1} \leq (\Delta^r)^{2r}$ elements of $D(w, r)$. For any element $A$ of $D(w, r)$, the probability of $f(A)$ differs from that of $A$ by a factor of at most $x^{2r}$. Hence, we have

(10) $$\mathbb{P}(S'(w, c, n, r)) \geq x^{2r} \Delta^{-2r^2} \mathbb{P}(D(w, r)).$$

Summing over $w$ and dividing by $|V_n|$ yields

(11) $$\mathbb{P}(S'(c, n, r)) \geq x^{2r} \Delta^{-2r^2} \mathbb{P}(D(r)).$$

However, the choice of $r$ (via Lemma 2.6) implies that if there are two large components, then at least $3/4$ of the vertices lie within distance $r$ of each one, so at least $1/2$ the vertices lie within distance $r$ of both. Hence, $\mathbb{P}(D(r)) \geq (1/2)\delta_n$. Combining this with (11) yields (8). □

PROOF OF THEOREM 2.1. By earlier remarks we may assume there is some $x > 0$ so that $x \leq p_n \leq 1 - x$ for all $n$. Now by Corollary 2.5, $\mathbb{P}(S'(c,n,r)) \to 0$ since $|V_n| \to \infty$. However, by Lemma 2.7, $\delta_n \leq 2x^{-2r}\Delta^{2r^2} \times \mathbb{P}(S'(c,n,r))$. Since $x$, $\Delta$ and $r$ are independent of $n$, it follows that $\delta_n \to 0$, precisely as we require. □

2.1. *Remarks.*

- For bond percolation on the complete graph with $n$ vertices, $G(n,p)$, it is known (see [8] or [17]) that for a suitable choice of $p$ (close to $1/n$) there are (roughly speaking) typically several components of order $n^{2/3}$; but whatever the choice of $p$, there is at most one component larger than this. It may be reasonable to strengthen Conjecture 1.1 in accordance with this.
- It turns out that we can strengthen Theorem 2.1 to give uniqueness of components of order $|V_n|^\omega$ for some $\omega < 1$. We need to allow the value of $r$ in the proof to vary with $n$, and we need to be considerably more careful in specifying how a configuration in $D(w,r)$, that is, one in which two large components intersect the ball $B(w,r)$, is transformed into a configuration containing an $L$-bridge. The details are given as Theorem 2.8.
- The value of $\omega$ implied by the proof of Theorem 2.8 can almost certainly be improved with more care; more difficult would be to establish the best possible value. Furthermore, one might expect rather more to be true, much as in the case of $G(n,p)$: roughly speaking, once the components become significantly larger than logarithmic in the number of vertices, they quickly agglomerate to form a single giant component. Therefore, except in a small window of values of $p_n$, one would expect at most one component of bigger than logarithmic size. In the case of $d$-regular expanders of high girth (see Section 3) we expect this window to be around $p_n = 1/(d-1)$. Note that we do know that for $p$ sufficiently close to 1 (independent of $n$), there is at most one component of greater than logarithmic size; see the remark after equation (13).
- The condition of expansion is a very strong one and it seems reasonable to conjecture that Theorem 2.1 holds under rather weaker conditions, such as some sublinear lower bound on the edge-boundary of subsets of the vertices. In fact, our proof does enable us to slightly weaken the expansion assumption, since the distance $r$ in the proof is allowed to grow (slowly) with $n$. In the context of vertex-transitive graphs, such a variant is Conjecture 1.1.
- See [22] for a proof using somewhat related ideas in a different context.

2.2. *A stronger uniqueness result.* In this section we show how to adapt our methods to establish the following stronger result.

THEOREM 2.8. *Given $b > 0$ and $\Delta \in \mathbb{N}$, there exists $\omega < 1$ such that the following holds for a sequence of vertex b-expanders $G_n = (V_n, E_n)$ with maximum degree at most $\Delta$, with $|V_n| \to \infty$. Let $0 \leq p_n \leq 1$. Then*

(12) $\quad \mathbb{P}(G_n(p_n) \text{ contains more than one component of order at least } |V_n|^\omega) \to 0$

*as $n \to \infty$.*

PROOF. We say that a component of $G_n(p_n)$ is *large* if it contains at least $|V_n|^\omega$ vertices; let $u_n = \lceil |V_n|^\omega \rceil$. We now imitate the arguments leading to Theorem 2.1 and see for what $\omega$ the proof still holds.

Lemma 2.2 is essentially unchanged: the probability that there is a component of size between $u_n$ and $|V_n|/2$ is bounded above by

$$\text{(13)} \quad \sum_{r=u_n}^{\lfloor |V_n|/2 \rfloor} \frac{|V_n|}{r} (\Delta e)^r (1-p_n)^{br} \leq \frac{|V_n|}{|V_n|^\omega} \frac{((\Delta e) A^b)^{|V_n|^\omega}}{1 - (\Delta e) A^b},$$

for $1 - A \leq p_n \leq 1$ much as in (3). This tends to 0 for any $\omega > 0$ (indeed, this even holds provided $u_n$ grows at least as fast as some particular multiple of $\log |V_n|$). Much as before, there are no large components for small $p$, so we can restrict to $p \in [x, 1-x]$.

The equivalent of Lemma 2.6 requires choosing $r_n$ so that if $|A| \geq u_n$, then $|B(A, r_n)| \geq (3/4)n$. Much the same argument as before shows that it is sufficient to choose

$$\text{(14)} \quad r_n = \lceil \log_{1+b}(|V_n|/u_n) \rceil = \lceil (1-\omega) \log_{1+b} |V_n| \rceil.$$

Now the definition of an $L$-bridge depends on the definition of a large component, which now depends on the choice of $\omega$. Much as in Corollary 2.4, we let $S(e, \omega, n)$ be the event that $e$ is an $L$-bridge in $G_n(p_n)$, and let $S(\omega, n)$ be the event that $S(e_n, \omega, n)$ occurs for an edge $e_n$ chosen uniformly at random. Then (6) becomes

$$\text{(15)} \quad \mathbb{P}(S(\omega, n)) \leq \frac{|V_n|}{|V_n|^\omega} \frac{\alpha}{\sqrt{|E_n|}}.$$

The proof of Corollary 2.5 is much unchanged, but since $r_n$ is no longer a constant the conclusion becomes

$$\mathbb{P}(S'(\omega, n, r_n)) \leq \frac{\Delta^{r_n+1}}{2} \mathbb{P}(S(\omega, n))$$

which, using (15), becomes

$$\text{(16)} \quad \leq \gamma \Delta^{r_n} |V_n|^{1/2-\omega},$$

for some $\gamma$ independent of $n$.

Greater care is needed in adapting Lemma 2.7. We need to reduce the $r^2$ appearing as an exponent in (8) to some multiple of $r_n$. In order to do this, we must be more precise about how we transform a configuration lying in $D(w, r_n)$ to one lying in $S'(w, \omega, n, r_n)$. Recall that a configuration is a subset, $F$ say, of the edge set $E_n$; we identify configurations with the corresponding spanning subgraph $(V_n, F)$, and a percolation process on $G_n$ is just a probability measure on the set of configurations. Recall also that the ball $B(w, r_n)$ is defined in terms of the original graph $(V_n, E_n)$.

Now, for each $w \in V_n$ and each unordered pair of vertices $x, y \in B(w, r_n)$, we fix one arbitrarily chosen path from $x$ to $y$, of length at most $2r_n$, lying entirely inside $D(w, r_n)$. Call this the *canonical path* $P(w, \{x, y\})$. Then, given a configuration, $F$, lying in $D(w, r_n)$, [i.e., a configuration such that at least two large components of $(V_n, F)$ intersect the ball $B(w, r_n)$] we obtain a configuration lying in $S'(w, \omega, n, r_n)$ as follows. Take two vertices, $x$ and $y$ say, that lie in $B(w, r_n)$, but which lie in different large components of $(V_n, F)$. Consider the process of adding, successively, the edges of the canonical path $P(w, \{x, y\})$ to the configuration $F$. (Note that some of these edges may already belong to $F$, and these are ignored in this process.) At some point the addition of one of these edges must join two large components. Stopping at this point (it does not matter whether before or after) makes this edge an $L$-bridge: this gives us our configuration lying in $S'(w, \omega, n, r_n)$.

As in the proof of Lemma 2.7, we have obtained a function, $f$ say, from $D(w, r_n)$ to $S'(w, \omega, n, r_n)$, but we have been more careful about the number of preimages each point can have. Since $B(w, r_n)$ contains at most $\Delta(\Delta - 1)^{r_n}/(\Delta - 2)$ vertices, it is not hard to see that it contains at most $5(\Delta - 1)^{2r_n}$ unordered pairs of vertices. Each such pair of vertices has a canonical path containing at most $2r_n$ edges and if $A \in D(w, r_n)$, then $A$ can be obtained from $f(A)$ by the deletion of some subset of the edges of some canonical path. Since a set of size no more than $2r_n$ has at most $2^{2r_n}$ subsets, we see that each configuration has at most $5(\Delta - 1)^{2r_n} 2^{2r_n}$ preimages. Just as in the proof of Lemma 2.7, the probability of $A$ differs from the probability of $f(A)$ by a factor of at most $x^{2r_n}$, so (10) becomes

$$\mathbb{P}(S'(w, \omega, n, r_n)) \geq \tfrac{1}{5} x^{2r_n} (\Delta - 1)^{-2r_n} 2^{-2r_n} \mathbb{P}(D(w, r)),$$

and, hence, (8) becomes

(17) $\qquad \mathbb{P}(S'(\omega, n, r_n)) \geq x^{2r_n} (\Delta - 1)^{-2r_n} 2^{-2r_n} \delta_n / 10.$

Combining (16) and (17) and simplifying a little gives

(18) $\qquad \delta_n \leq 10\gamma \Delta^{3r_n} x^{-2r_n} 2^{2r_n} |V_n|^{1/2 - \omega}.$

Recalling the choice of $r_n$, (14), we obtain

$$\delta_n \leq 10\gamma (4\Delta^3 x^{-2})^{(1-\omega)\log_{1+b} |V_n| + 1} |V_n|^{1/2 - \omega}$$
$$= 10\gamma (4\Delta^3 x^{-2}) |V_n|^{(1-\omega)\log_{1+b}(4\Delta^3 x^{-2})} |V_n|^{1/2 - \omega}.$$

So we see that $\delta_n \to 0$ as $n \to \infty$ provided

$$(19) \qquad (1-\omega)\log_{1+b}(4\Delta^3 x^{-2}) + (\tfrac{1}{2} - \omega) < 0.$$

Since (19) clearly holds for $\omega$ sufficiently close to 1, this establishes that $\delta_n$—the probability of two or more large components—tends to 0 as $n \to \infty$ for such values of $\omega$, exactly as we require. $\square$

**3. High girth expanders.** In this section we show that when we consider $d$-regular expanders of girth tending to infinity, we can identify the critical value of $p$ above which a (unique) giant component appears, namely, $1/(d-1)$.

PROPOSITION 3.1. *For every $\varepsilon > 0$, there is an $a = a(\varepsilon) > 0$ and a $\delta = \delta(\varepsilon) > 0$ such that the following holds. Let $G = (V, E)$ be a finite, $d$-regular graph on a set $V$ of $n$ vertices, let $g = g(G)$ denote its girth and let $c = c(G)$ denote its isoperimetric number. If*

$$C = \frac{c}{2}\left(\frac{\varepsilon}{2d}\right)^{d/ca} - \frac{3\ln 2}{(1+\varepsilon/3)^{g/2}} > 0,$$

*then, for $p = \frac{1+\varepsilon}{d-1}$, the random graph $G(p)$ has, with probability that exceeds $1 - e^{-Can} - e^{-\delta n}$, a connected component with at least $an$ vertices.*

PROOF. It is convenient to consider the random subgraph $G(p)$ as a union of two independently chosen random subgraphs $G(p_1)$ and $G(p_2)$, where $p_1 = \frac{1+\varepsilon/2}{d-1}$ and $p_2$ ($\geq \frac{\varepsilon}{2d}$) is chosen such that $(1-p_1)(1-p_2) = 1-p$. Seen from any vertex, the graph $G$ out to a distance $g/2$ looks just like a $d$-regular tree. By standard results from the theory of branching processes (see, e.g., [16]), with probability at least $1 - e^{-\delta(\varepsilon)n}$, at least $a' = a'(\varepsilon)n$ vertices of $G(p_1)$ lie in components of size at least $m$, where $m = (1+\varepsilon/3)^{g/2}$. Conditional on this, define $a = a'/3$ and fix a set of at most $\frac{a'n}{m}$ such components that contain together at least $a'n$ vertices. We claim that with probability at least $1 - e^{-Can}$, in the random graph $G(p_2)$ there is no way to split these components into two parts $A$ and $B$, each containing at least $a'n/3$ vertices, with no path of $G(p_2)$ connecting the two parts. This will imply that with the required probability, the union of the two graphs $G(p_1)$ and $G(p_2)$ contains a connected component consisting of at least $a'n/3 = an$ vertices, as needed.

To prove the claim notice first that there are at most $2^{a'n/m}$ possible ways to split the components into two sets $A$ and $B$ as required. For each such fixed choice, the fact that $c = c(G)$ and Menger's theorem imply that there are at least $ca'n/3$ pairwise edge-disjoint paths in $G$ from $A$ to $B$. As $G$ has $dn/2$ edges, at least half of these paths are of length at most $3d/(ca')$ each. The probability that none of those paths belongs to $G(p_2)$ is at most

$$(1-p_2^{3d/(ca')})^{ca'n/6} \leq \left[1-\left(\frac{\varepsilon}{2d}\right)^{3d/(ca')}\right]^{ca'n/6} \leq \exp\left(-\frac{ca'n}{6}\left(\frac{\varepsilon}{2d}\right)^{3d/(ca')}\right).$$

It follows that the probability that there is some partition into sets $A$ and $B$ as above, with no paths of $G(p_2)$ between them, is at most

$$2^{a'n/m}\exp\left(-\frac{ca'n}{6}\left(\frac{\varepsilon}{2d}\right)^{3d/(ca')}\right) = e^{-Can},$$

which completes the proof. □

Simple branching process comparisons show that on any $d$-regular graph, $G$, if $p < 1/(d-1)$, then the probability that $G(p)$ has a large component is small. Combining this fact with the above proposition easily gives the following theorem, which can be loosely described as saying that the critical probability for the emergence of a giant component, in a sequence of $d$-regular expanders with girth tending to infinity, is $1/(d-1)$.

THEOREM 3.2. *Let $d \geq 2$ and let $(G_n)$ be a sequence of $d$-regular expanders with girth $(G_n) \to \infty$.*

*If $p > 1/(d-1)$, then there exists $c > 0$ such that:*

$\mathbb{P}(G_n(p) \text{ contains a component of order at least } c|V(G_n)|) \to 1 \qquad \text{as } n \to \infty.$

*If $p < 1/(d-1)$, then for any $c > 0$,*

$\mathbb{P}(G_n(p) \text{ contains a component of order at least } c|V(G_n)|) \to 0 \qquad \text{as } n \to \infty.$

REMARKS.

- The arguments above imply that for every fixed $d$, the critical probability for the emergence of a linear size connected component in a *random d-regular graph* on $n$ vertices is almost surely $1/(d-1) + o(1)$. Indeed, these graphs do have some constant size cycles, but their number is, almost surely, small enough that they can be ignored, including near critical behavior.
- It might be possible to apply the techniques above and show that vertex transitive graphs of degree $d$ in which the girth is proportional to the diameter also have the same critical probability, since it is known (see, e.g., [5]) that such graphs are good expanders. A rigorous proof may require some care, as the proposition itself does not suffice here. On the other hand, there are simple examples showing that without the assumption of vertex transitivity, the conclusion fails. A counterexample can be constructed by taking some $\log n$ 3-regular expanders, each on $n$ vertices and each of logarithmic girth, by omitting an edge from each of them, and then by joining them all along a cycle keeping the resulting graph 3-regular.
- One can use the approach above to prove that the giant component in the Erdős–Rényi random graph $G(n,p)$ emerges at $p = 1/n$.

As we remarked in Section 2, the condition of expansion is rather strong, and one would expect similar results to hold under weaker conditions. In the context of transitive graphs, we suggest the following conjecture.

CONJECTURE 3.3. Let $\{G_n\}_{n \in N}$ be a sequence of $d$-regular connected finite transitive graphs, $|V_n| \nearrow \infty$, and suppose that $\mathrm{diameter}(G_n) = o(|V_n|/\log|V_n|)$. Then the threshold for the existence of a connected component of size $|V_n|/10$, with probability $1/2$, is uniformly bounded away from 1.

By [12] the threshold is sharp. Note that the conjecture is true for tori; see, for example, the section on percolation in a wedge in [14]. Recently the conjecture has been shown to hold for certain Cayley graphs [23]. However, for general graphs, even if we make a stronger assumption that

$$\mathrm{diameter}(G_n) < |V_n|^\varepsilon,$$

for some $\varepsilon < 1$, we do not know how to show that the threshold for a giant component is bounded away from 1.

**4. Percolation in $\mathbb{Z}^d$.** In 1990 Kesten [18] proved that the critical probability for bond percolation in $\mathbb{Z}^d$ is $\frac{1+o(1)}{2d}$, where the $o(1)$ error term tends to zero as $d$ tends to infinity. Hara and Slade [15] obtained a better estimate for the error term, and Bollobás and Kohayakawa [10] gave a somewhat simpler proof. Here we sketch a simpler argument giving the result of Kesten, following the method of [1].

The fact that $1/(2d-1)$ is a lower bound for the critical probability is obvious, hence, we only sketch the proof of the upper bound. It will be convenient to prove the upper bound for the subgraph $G$ of $\mathbb{Z}^d$ induced on $\mathbb{Z}^2 \times [d]^{d-2}$, where $[d] = \{1, 2, \ldots, d\}$. We assume, from now on, that $d$ is sufficiently large. Let $\varepsilon > 0$ be small, and put $p = \frac{1+\varepsilon}{2d}$. It is convenient to first consider, in phase one, the random subgraph $G(p)$ of $G$ obtained by taking each edge, randomly and independently, with probability $p$, and take, in phases two and three, its union with two additional randomly chosen subgraphs $G(p_i)$ with, say, $p_i = \frac{1}{d^2}$ for each of them.

Split the vertex set of $G$ into $d$-dimensional boxes, each isomorphic to $[d]^d$. Each two neighboring boxes have $d^{d-1}$ edges connecting them. The basic approach is to show that in each fixed box, after the first two phases our random subgraph will have, with high probability, a linear size component with lots of neighbors in the boundary. The result can then be obtained by taking the additional fresh random edges of phase three and by using some very rough estimates on percolation in $\mathbb{Z}^2$.

We first need some (known) expansion properties of $[d]^d$. Since what we need is extremely simple, we include a proof.

LEMMA 4.1. *The edge-isoperimetric number of the graph $[d]^d$ is at least $1/(2d)$. That is, for every set $A$ of at most half the vertices of $[d]^d$, there are at least $\frac{|A|}{2d}$ edges connecting $A$ to its complement.*

PROOF. For every pair of vertices $a \in A$ and $b \notin A$, take a canonical path from $a$ to $b$ obtained by changing the coordinates in which $a$ and $b$ differ one by one, from left to right, where each coordinate is being changed monotonically. Each such path must contain an edge connecting a vertex of $A$ with one in its complement, and every edge appears in at most $d^{d+1}$ paths. Therefore, there are at least $|A|(d^d - |A|)/d^{d+1}$ edges connecting $A$ to its complement. □

Consider, now, a random subgraph $H(p)$ with $p$ as above, where $H$ is the subgraph of $G$ induced by $[d]^d$. Let $\delta > 0$ be a fixed small real (smaller than $\varepsilon/2$, say). Call a vertex of $H$ *good* if it has at most $\delta d/10$ coordinates which are either 1 or $d$. Note that each such vertex has at least $(2 - 2\delta/10)d$ neighbors inside $H$. Put $n = d^d$. Call a connected component of $H(p)$ an *atom* if it has at least, say, $d^{100}$ vertices. We first claim that with high probability, by the end of phase one, every vertex of $H(p)$, besides at most some $n/2^{c_1 d}$, has at least one neighbor which lies in $H(p)$ in an atom, where here $c_1 = c_1(\delta) > 0$. Indeed, all vertices but some $n/2^{c_2 d}$ are good. Each such vertex has at least $(1 - \frac{\delta}{10})d$ coordinates that are neither 1 nor $d$. Without loss of generality assume these are the first coordinates. Let $v = (v_1, \ldots, v_d)$ be the vertex. For each $i \leq \delta d/5$ (say), consider the connected component of the neighbor $(v_1, v_2, \ldots, v_{i-1}, v_i + 1, v_{i+1}, \ldots, v_d)$ of $v$ obtained by considering the (forward) branching process only on vertices of the form $(v_1, v_2, \ldots, v_i + 1, u_{i+1}, u_{i+2}, \ldots)$ with each $u_j$ for $j > i$ being in the set $\{v_j, v_j - 1, v_j + 1\}$. In this process we always move from a vertex to ones with bigger Hamming distance from $v$, and if a vertex is obtained more than once as a child, we omit it. This is done some $c(\delta) \log d$ generations. An easy calculation shows that with probability tending to 1 (as $d \to \infty$), there are no vertices obtained more than once as a child (and thereby omitted); then, standard results on branching processes imply that for each fixed neighbor, we manage to grow an atom with probability bounded away from zero. As the events for distinct neighbors we are considering are independent, the desired claim follows.

Consider, now, the set of all atoms obtained. These cover together a constant fraction, say $c_3 n$, of the $n$ vertices of $H$ [with $c_3 = c_3(\delta)$]. Now add, in phase two, edges of $H$ randomly, with probability $\frac{1}{d^2}$. We claim that in the resulting graph, with high probability, no union of atoms $A$ covering at least, say, $n/d^5$ vertices can be separated from the union $B$ of all other atoms, when this union also covers at least $n/d^5$ vertices. To prove this claim, denote, for any set $X$ of vertices of $H$, by $N(X)$ the set of all its neighbors in $H$. Consider two possible cases.

CASE 1. $|N(A) \cap N(B)| \geq \frac{n}{d^{10}}$. In this case there are at least $\frac{n}{d^{10}}$ pairwise edge disjoint paths of length 2 connecting $A$ and $B$, and the probability none of them is chosen is at most

$$\left(1 - \frac{1}{d^4}\right)^{n/d^{10}}.$$

Even when multiplied by the number of possibilities for choosing $A$ and $B$, which is smaller than $2^{n/d^{100}}$, this is negligible.

CASE 2. $|N(A) \cap N(B)| < \frac{n}{d^{10}}$. Assume, without loss of generality, that $|B| \geq |A|$. Since $A \cup N(A)$ misses most of $B$, it is not very large, and by Lemma 4.1 we get that there are at least some $c_4 n/d^7$ distinct vertices of $H$ of distance 2 (in $H$) from $A$ [because there are at least $c_4 n/d^6$ edges connecting $A \cup N(A)$ to its complement]. Most of these vertices have neighbors that are atoms and, hence, lie in $B$. This gives many paths of length 3 between $A$ and $B$, and it is easily seen this construction gives $\Omega(n/d^8)$ pairwise edge disjoint paths from $A$ to $B$. As before, with high probability, all the edges of at least one of those will be chosen in phase two.

The preceding argument establishes the claim. Moreover, it implies that with high probability, by the end of the second phase there is a connected component of the resulting graph that contains all the vertices that were in atoms by the end of the first phase, besides at most $n/d^5$ of them. Let us call this component the *distinguished* component. Note also that with high probability, all vertices of $H$, besides at most some $O(n/d^4)$, have at least one neighbor in this component.

We can now consider each copy of $H$ among the ones that split the vertices of $G$ as a site which, with probability very close to 1, is present; where we make it *present* if on each face of its boundary at least 3/4 of the vertices have neighbors in the distinguished component. By the above discussion this happens with high probability, and each copy of $H$ behaves independently in this respect. Taking now, in phase three, fresh edges with probability $\frac{1}{d^2}$, we get that with extremely high probability, every two such neighboring sites become connected, as there are at least $d^{d-1}/2$ potential pairwise disjoint paths of length 2 that connect the two corresponding distinguished components, and with (very) high probability, at least one of those will be chosen in the third phase. Call a copy of $H$, regarded as a site in its own right, *active* if it is present and it is connected to all its present neighbors in the way just described; so for $d$ sufficiently large the probability that a site is active is arbitrarily close to 1. Sets of sites at ($l_\infty$) distance greater than 1 behave independently, and it is well known that 1-dependent bond percolation on $\mathbb{Z}^2$ with a sufficiently high marginal percolates (for much stronger results see [19]); hence, for $d$ sufficiently large the active sites connect up to form an infinite component.

**5. Concluding remarks.** It seems plausible that if $G = (V, E)$ is an expander with $n$ vertices, and $p$ is above the critical probability for the emergence of a giant component, then the giant component of $G(p)$ will have, itself, reasonably strong expansion properties. The proof of Lemma 2.2 can be easily modified to prove that this is, indeed, the case at least when $p$ is close to 1.

PROPOSITION 5.1. *Let $b > 0$ and let $G_n = (V_n, E_n)$ be a sequence of graphs with maximum degree at most $\Delta$ which are edge $b$-expanders, with $|V_n| \to \infty$. Let $A > 0$ be such that $(\Delta e) 2^b A^{b/2} < 1/2$ and let $1 - A \le p_n \le 1$ for each $n$. Then*

$$
(20) \qquad \mathbb{P}\bigg( G_n(p_n) \text{ is not a } \frac{1}{\log_2 n} \text{ edge expander} \bigg) \to 0
$$

*as $n \to \infty$.*

SKETCH OF PROOF. A simple modification of the proof of Lemma 2.2 shows that the probability that there is a connected induced subgraph of $G_n(p_n)$, whose size $r$ is at least $\log_2 n$ and at most $n/2$, which has at most $br/2$ edges emanating from it to the rest of the graph is at most

$$
\sum_{r=\log_2 n}^{n/2} \frac{n}{r} (\Delta e)^r \binom{br}{br/2} (1-p_n)^{br/2} \le \frac{n}{\log_2 n} \sum_{r \ge \log_2 n} (\Delta e 2^b A^{b/2})^r < \frac{2}{2\log_2 n}.
$$

The desired result follows. □

Consider an infinite transitive graph $G$. It is believed that uniqueness of the infinite cluster holds at all $p$ iff $G$ is amenable, see [7]. For the product of an infinite regular tree and $\mathbb{Z}$, and certain other nonamenable graphs, nonuniqueness of the infinite cluster at some values of $p$ is known to hold. We suspect that on finite transitive graphs this is not the case in the following sense.

CONJECTURE 5.2. *Let $\{G_n\}_{n \in N}$ be a sequence of $d$-regular connected finite transitive graphs, $|V_n| \nearrow \infty$. Given $\delta > 0$, there is $C > 0$ so that if for any $n$ and $v \in G_n$,*

$$\mathbb{P}_p(v \text{ is in a component of diameter} \geq \mathrm{diam}(G_n)/2) > \delta,$$

*then*

$$\lim_n \mathbb{P}_p(\text{there is a connected component of size} \geq C|V_n|) = 1.$$

That is, once a fixed vertex is with positive probability in a large cluster in the sense of diameter, it will be in a cluster which is large in the sense of volume. Note that the conjecture is true for Euclidean lattice tori and, by the discussion above, for expanders with growing girth. A useful fact that supports the conjecture is that finite transitive graphs do not admit bottle necks. The Cheeger constant of a finite transitive graph is at least the reciprocal of the diameter (see [5]).

**Acknowledgment.** Part of this work was carried out during a visit of the first and last authors at Microsoft Research, Redmond, Washington. The first and last authors would like to thank their hosts at Microsoft for their hospitality. We also thank Ehud Freidgut for helpful comments.

SCHOOLS OF MATHEMATICS
AND COMPUTER SCIENCE
RAYMOND AND BEVERLY SACKLER
FACULTY OF EXACT SCIENCES
TEL AVIV UNIVERSITY
TEL AVIV 69978
ISRAEL
E-MAIL: noga@math.tau.ac.il

MICROSOFT RESEARCH
AND WEIZMANN INSTITUTE
ONE MICROSOFT WAY
REDMOND, WASHINGTON 98052
USA
E-MAIL: itai@wisdom.weizmann.ac.il

DEPARTMENT OF PURE MATHEMATICS
AND MATHEMATICAL STATISTICS
CENTRE FOR MATHEMATICAL SCIENCES
WILBERFORCE ROAD
CAMBRIDGE
UNITED KINGDOM
E-MAIL: AM.Stacey@dpmms.cam.ac.uk